\documentclass[final,1p,times,english]{article}
\usepackage{color,latexsym,amsmath,amsfonts,amssymb}
\usepackage[utf8]{inputenc}
\usepackage{layout,verbatim,amssymb,amsmath,setspace}
\usepackage{amsthm}

\newcommand{\vG}{\varGamma}
\newcommand{\ve}{\varepsilon}
\newcommand{\wt}{\widetilde}

\newcommand{\ov}{\overline}
\newcommand{\mg}{\marginpar}

\newcommand{\vp}{\varepsilon}
\def\N{{\mathbb N}}

\def\mcL{{\mathcal L}}
\def\mcS{{\mathcal S}}
\def\mcI{{\mathcal I}}

\def\mcP{{\mathcal P}}

\def\emp{\emptyset}

\def\be{\begin{equation}}
\def\ee{\end{equation}}
\def\ba*{\begin{eqnarray*}}
\def\ea*{\end{eqnarray*}}

\def\bs{\boldsymbol}
\def\vPh{\varPhi}

\newtheorem{thm}{Theorem}[section]
\newtheorem{deff}[thm]{Definition}

\newtheorem{rem}[thm]{Remark}
\newtheorem{prop}[thm]{Proposition}
\newtheorem{cor}[thm]{Corollary}
\newtheorem{que}[thm]{Question}
\newtheorem {ex}[thm]{Example}

\begin{document}
\hyphenation{pro-per-ties} \hyphenation{cha-rac-te-ri-za-tions}
\hyphenation{stron-gly} \hyphenation{pro-per-ties}
\hyphenation{Pro-per-ty} \hyphenation{si-mi-lar}
\hyphenation{re-gu-la-ri-ty} \hyphenation{ine-qua-li-ty}
\hyphenation{non-over-lapping} \hyphenation{glo-bal-ly}
\hyphenation{sub-in-ter-vals} \hyphenation{se-quen-ce}
\hyphenation{u-ni-for-mly} \hyphenation{ge-ne-ra-li-zed}
\hyphenation{con-ti-nu-ous} \hyphenation{glo-bal-ly}
\hyphenation{sub-in-ter-vals} \hyphenation{se-quen-ce}
\hyphenation{u-ni-for-mly} \hyphenation{ge-ne-ra-li-zed}
\hyphenation{ter-mi-no-lo-gy} \hyphenation{re-fe-ren-ce}
\hyphenation{Theo-rem} \hyphenation{sub-in-ter-vals}
\hyphenation{dif-fe-ren-tia-ble} \hyphenation{u-ni-for-mly}
\hyphenation{e-qui-va-lent} \hyphenation{no-ti-cing}
\hyphenation{ge-ne-ra-ted}

\begin{center}
 {\bf \Large  Gauge integrals and selections\\ of  weakly compact valued multifunctions
\footnote{D. Candeloro, A.R. Sambucini: Department of Mathematics and Computer Sciences - 06123 Perugia (Italy), email: domenico.candeloro@unipg.it, anna.sambucini@unipg.it\\
Di Piazza: Department of Mathematics, University of Palermo,  Via Archirafi 34, 90123 Palermo (Italy), email: luisa.dipiazza@unipa.it\\
Musial: Institut of Mathematics, Wroc{\l}aw University,  Pl. Grunwaldzki  2/4, 50-384 Wroc{\l}aw (Poland), email: musial@math.uni.wroc.pl
}
}
\\
\vspace{0.2cm} D. Candeloro,  L. Di Piazza, K. Musia{\l} and A. R. Sambucini\footnote{orcid id:   D. Candeloro: 0000-0003-0526-5334, 
A.R. Sambucini: 
0000-0003-0161-8729} 
\end {center}
\title{\bf}
\author{ D. Candeloro, L. Di Piazza, K. Musia{\l} and A. R. Sambucini}
\begin{abstract}
In the paper  Henstock, McShane, Birkhoff and variationally multivalued integrals are studied
for multifunctions taking values in the hyperspace of convex and weakly compact subsets of a general Banach space $X$. In particular the existence of selections integrable in the same sense of the corresponding multifunctions has been considered.
\end{abstract}

\noindent \\
\begin{msc}
  28B20,  26E25, 26A39,
 28B05,  46G10, 54C60, 54C65
\end{msc}

\section{Introduction}
After the pioneering works of Aumann  and Debreu, the theory of multivalued functions  has been intensively studied and  several notions of integral have been developed using different techniques, see for example \cite{ph,ms2001,ms2002,CASCALES2,bs2004,dm7,ckr1,ckr,ckr2,dm2,mu4,BS2011,bms,dm5,dm,cgsub1,BCS,ccgs2015}.
These notions have shown to be useful when modeling some theories in different fields as  optimal
control and mathematical economics, see for example \cite{hp,ebs2005,cc,dps1,dms}.
The choice to deal with these types of integration is motivated
 by the fact that a very important tool in this framework is the Kuratowski
and Ryll-Nardzewski
theorem which guarantees the existence of measurable selectors, though this famous result   requires the separability of the range space $X$.

The starting point of this research are  the papers   \cite{ckr1,ckr,ckr2,mu4} in which this result
 was extended to the non separable Banach spaces  for the Pettis multivalued integral and the papers \cite{dm,dm2,dm5,dm7} where the existence of selections integrable in the same sense of the corresponding multifunctions has been considered for some gauge integrals  in the hyperspace  $cwk(X)$  ($ck(X)$) of convex and weakly compact (compact) subsets  of a general Banach space $X$
and \cite{BCS,cs2014,cs2015} in which these arguments are studied in suitable Banach lattices.

In the present paper selection results are obtained for Henstock, McShane, Birkhoff and variational integrals of  $cwk(X)$-valued multifunctions. Also  a number of examples are given when this is not possible.

The paper is organized as follows: in section \ref{two}  the Henstock, McShane, Birkhoff and variationally multivalued integrals are considered
for multifunctions taking values in  $cwk(X)$, and an embedding result, which will be useful later, is recalled. Moreover some other properties of variational integrals are proved,  as for example their Bochner measurability.

Section \ref{five} is more specifically devoted to selections theorems: results concerning McShane (resp. Henstock) integrability of all selectors of a McShane (resp. Henstock) integrable multifunction are given, at least in the separable case, and examples have been given to show that in non-separable spaces the results can fail.

In case of a variationally Henstock (resp. McShane) integrable $cwk(X)$-valued multifunction, the problem of existence of at least one variationally Henstock (resp. McShane) selection has been investigated. Full solution has been given for Banach spaces with  the Radon-Nikod\'ym Property. In case of  a general space $X$ we present a solution for the variational McShane integral. It remains an open problem if, in a general Banach space $X$, any $cwk(X)$-valued variationally Henstock-integrable multifunction has a variationally Henstock-integrable selection. Results with positive answers are given in some particular cases.

In section \ref{three}, thanks to the structure of near vector space of $cwk(X)$, it is shown  that the variationally Henstock integrability of $F$, with \mbox{$0 \in F$}, implies its Birkhoff  integrability and this result is used in order to obtain a
 decomposition  extending the one given in  \cite{dm} for compact convex valued multifunctions.
For scalarly measurable multifunctions some implications are also given, connecting
 Henstock, McShane and Pettis multivalued integrability to one another. Moreover an example of a variationally Henstock but not variationally McShane integrable multifunction is given.

\section{Preliminary facts}\label{two}

Throughout $[0,1]$ is the unit interval of the real line equipped
with the usual topology and  Lebesgue measure $\lambda$,
$\mathcal{L}$ denotes the family of all Lebesgue measurable subsets of $[0,1]$, and
 $\mcI$
is the collection of all closed subintervals of $[0,1]$.  If $I\in \mcI$
then $|I|$ denotes its length.

$X$ is an arbitrary  Banach space with its dual $X^*$.
  The closed unit ball of $X^*$ is denoted by
  $B(X^*)$. $cwk(X)$ is    the family of all non-empty  convex weakly compact subsets of $X$ and  $ck(X)$ is the family of all  compact members of $cwk(X)$ .  We consider on $cwk(X)$ the
   Minkowski addition ($A+B :\,=\{a+b:a\in A,\,b\in B\}$) and the
   standard multiplication by scalars. $\|A\|:=\sup\{\|x\|\colon x\in{A}\}$.
$d_H$ is the Hausdorff distance in $cwk(X)$. $cwk(X)$ with the Hausdorff
   distance is a complete metric space.
\\
     For every $C \in cwk(X)$ the {\it support
  function of} $C$ is denoted by $s( \cdot, C)$ and defined on $X^*$ by $s(x^*, C) = \sup \{ \langle x^*,x \rangle : \ x
  \in C\}$, for each $x^* \in X^*$. For all unexplained definitions we refer to \cite{CV}.\\

  A map $\vG:[0,1]\to 2^X\setminus\{\emptyset\}$
  (= non-empty subsets of $X$) is called a {\it multifunction}.
 A multifunction $\vG:[0,1]\to 2^X\setminus\{\emptyset\}$ is said to be a {\it simple multifunction} if there exists a finite collection $\{E_1,..., E_p\}$ of measurable subsets of $[0,1]$, pairwise disjoint, such that $\vG$ is constant on each $E_j$.
\\
 A multifunction $\vG:[0,1]\to cwk(X)$ is said to be {\it Effros measurable} (or simply {\it measurable})  if for each open $O \subset X$, the set $ \{t\in [0,1]: \vG(t)\bigcap O \neq \emptyset\}$ is measurable.
\\
 A multifunction $\vG:[0,1]\to cwk(X)$ is said to be {\it scalarly measurable} if for every $ x^* \in X^*$, the map
  $s(x^*,\vG(\cdot))$ is measurable. $\vG$ is said to be {\it Bochner measurable}  if there exists a sequence of simple multifunctions $\vG_n: [0,1] \to cwk(X)$ such that
$$\lim_{n\rightarrow \infty}d_H(\vG_n(t),\vG(t))=0$$ for almost all $t \in [0,1]$.\\

It is well known that the measurability of a $cwk(X)$-valued multifunction yields its scalar measurability (if $X$ is separable also the reverse implication is true). Moreover,
each Bochner measurable multifunction is also measurable (see \cite{hp}). The reverse implication fails (see Example \ref{ex10}).\\

A measurable multifunction $\vG$ is said to be {\it Aumann integrable} if it admits at least one Bochner integrable selection.
  A function $f:[0,1]\to X$ is called a {\it selection of} $\vG$ if $f(t) \in\vG(t)$,   for every $t\in [0,1]$.\\

 A {\it partition} ${\mathcal P}$ {\it in} $[0,1]$ is a collection $\{(I_1,t_1),$ $ \dots,(I_p,t_p) \}$,
  where $I_1,\dots,I_p$ are nonoverlapping subintervals of $[0,1]$, $t_i$ is a point of $[0,1]$, $i=1,\dots,p$.
 If $\cup^p_{i=1}I_i=[0,1]$, then  ${\mathcal P}$ is {\it a partition of} $[0,1]$. If   $t_i$ is a point of $I_i$, $i=1,\dots,p$,  we say that $P$ is a {\it Perron partition of}
  $[0,1]$.\\
 A {\it gauge} on $[0,1]$ is a positive function on $[0,1]$. For a given gauge $\delta$ on $[0,1]$,
  we say that a partition $\{(I_1,t_1), \dots,(I_p,t_p)\}$ is $\delta$-{\it fine} if
  $I_i\subset(t_i-\delta(t_i),t_i+\delta(t_i))$, $i=1,\dots,p$.
\\
 Given a function $g\colon [0,1] \to X$ and a partition $ {\mathcal P}=\{(I_1,t_1),$ $ \dots,(I_p,t_p) \}$  in $[0,1]$ we set
$$ \sigma(g,{\mathcal P})= \sum_{i=1}^p|I_i|g(t_i).$$

    \begin{deff} \rm
 A multifunction $\vG:[0,1]\to cwk(X)$ is said to be {\it Henstock} (resp. {\it McShane})
   integrable on $[0,1]$,  if there exists a non empty closed convex set  $\vPh_{\vG}([0,1])\subset{X}$
    such that for every $\varepsilon > 0$ there exists a gauge $\delta$ on $[0,1]$
such that for each $\delta$--fine Perron partition (resp. partition)
   $\{(I_1,t_1), \dots,(I_p,t_p)\}$ of
   $[0,1]$, we have
\begin{eqnarray}\label{e14}
d_H \left(\vPh_{\vG}([0,1]),\sum_{i=1}^p\vG(t_i)|I_i|\right)<\ve\,.
\end{eqnarray}\\
A multifunction $\vG:[0,1]\to cwk(X)$ is said to be {\it Birkhoff} (resp. {\it abs-Birkhoff})
   integrable on $[0,1]$,
   if there exists a non empty closed convex set  $\vPh_{\vG}([0,1]) \in cwk(X)$
 with the following property: for every $\varepsilon > 0$  there is a countable
partition $\Pi_0$  of $[0,1]$  in $\mathcal{L}$ such that for every countable partition $\Pi = (A_n)_n$
 of $[0,1]$  in $\mathcal{L}$
finer than $\Pi_0$ and any choice $T = (t_n)_n$  in $A_n$, the series
$\sum_n |A_n| \vG(t_n)$
 is unconditionally convergent  (resp. absolutely convergent)  and
\begin{eqnarray}\label{e14-a}
d_H \left(\vPh_{\vG}([0,1]),\sum_n \vG(t_n)|I_n|\right)<\ve\,.
\end{eqnarray}
\\
A multifunction $\vG:[0,1]\to cwk(X)$ is said to be {\it Henstock} (resp. {\it McShane or Birkhoff})
   integrable on $I\in\mcI$ if $\vG 1_I$ is respectively integrable on $[0,1]$. We write then $(H)\int_I\vG\,dt:=\vPh_{\vG 1_I}([0,1])$ (resp.  $(MS)\int_I\vG\,dt:=\vPh_{\vG 1_I}([0,1]$) or $(Bi)\int_I\vG\,dt:=\vPh_{\vG 1_I}([0,1]$)).
 It is known that a multifunction that is Henstock (McShane, Birkhoff) integrable on $[0,1]$ is in the same manner integrable on each $I\in\mcI$ (see e.g. \cite{dm}).
\end{deff}
  It is easily seen from the definition and the completeness of the Hausdorff metric that $ck(X)$ ($cwk(X)$)-valued  integrable multifunctions have compact (weakly compact) values of their integrals.

Moreover we would like to recall that
a vector function $g:[0,1]\to
X$ is Birkhoff-integrable if it is McShane integrable,
but just measurable gauges are involved in the notion of McShane integrability.
More precisely
\begin{deff}\rm
 $g:[0,1] \to X$ is Birkhoff-integrable if and only if
there exists an element $y\in X$ such that
for each $\vp>0$ a {\em measurable} gauge $\delta$ can be found on $[0,1]$, such that,
as soon as ${\mathcal P}=\{(t_j,I_j):j=1,...,n\}$ is any $\delta$-fine
partition of $[0,1]$, it holds
$\| \sigma(g,{\mathcal P}) - y \| \leq \vp$.
\end{deff}

(For the equivalence of this definition with the more
common notion of Birkhoff integrability see \cite{
nara,
ccgs2015})

\begin{deff}\rm

   A multifunction $\vG:[0,1]\to cwk(X)$ is said to be {\it variationally Henstock} ({\it variationally McShane})
   integrable,
   if there exists a finitely additive multifunction  $\vPh_{\vG}: {\mathcal I} \to {cwk(X)}$ with the following property:
   for every $\ve>0$ there exists a gauge $\delta$
   on $[0,1]$ such that for each $\delta$--fine Perron partition (partition)
   $\{(I_1,t_1), \dots,(I_p,t_p)\}$ of \/
   $[0,1]$, we have
\begin{eqnarray}\label{aa}
\sum_{j=1}^pd_H \left(\vPh_{\vG}(I_j),\vG(t_j)|I_j|\right)<\ve\,.
\end{eqnarray}
   We write then $({\rm v}H)\int_0^1\vG\,dt:=\vPh_{\vG}([0,1])$ ($(vMS)\int_0^1\vG\,dt:=\vPh_{\vG}([0,1]$)).  We call the set  multifunction  $\vPh_{\vG}$ the {\it variational Henstock} ({\it McShane}) {\it primitive} of $\vG$.
The variational integrals on $I\in\mcI$ are defined in the same way
as the ordinary ones. The integrals are uniquely determined.
\\
We say that a multifunction  $\vG$ is  {\it scalarly Henstock}
   {\it integrable} if, for every $x^* \in X^*$, the function $s(x^*,\vG(\cdot))$ is {\it Henstock
   integrable}.
   \end{deff}

 It follows from the definitions  that if $\vG$ is  McShane (variationally McShane) integrable, then it is also
    Henstock (variationally Henstock) integrable (with the same values of the integrals).   There is a McShane integrable function $f:[0,1]\to{B(L_{\infty}[0,1])}$ that is  not variationally McShane integrable \cite[3D]{fm}.
 \\

   When a multifunction is a function $f: [0,1]\to X$, then the set $\vPh_f([0,1])$ is reduced to a vector of $X$ and the above definitions coincide with those corresponding  for vector valued functions.
 For vector valued functions on $[0,1]$ the variational McShane integrability is equivalent to the Bochner
    integrability \cite{cx}.  This result has been generalized in \cite{dm11} to the case of vector valued functions defined in a compact finite Radon measure space.\\

 For the definitions of Pettis  and  of Henstock-Kurzweil-Pettis  integral for multifunctions we refer the reader  to
\cite{mu,mu3,mu4,mu8,ckr1,dm1,dm2}.\\

   In \cite{L1} a R{\aa}dstr\"{o}m embedding theorem  is extended to the space $cwk(X)$  as follows.

\begin{thm}
{\rm  (\cite[Theorem 5.6]{L1})}
 \label{5.6}
There exist a compact Hausdorff space $\Omega$ and a  map
$i:cwk(X) \to C(\Omega)$ such that
\begin{description}
\item[\ref{5.6}.1)] $i(\alpha A+ \beta C) = \alpha i(A) + \beta i(C)$ for every $A,C\in cwk(X), \alpha, \beta \in
\mathbb{R}^+;$
\item[\ref{5.6}.2)] $d_H(A,C)=\|i(A)-i(C)\|_{\infty},\ \ A,C\in cwk(X)$;
\item[\ref{5.6}.3)] $i(cwk(X))=\ov{i(cwk(X))}\ \ \ \ {\rm (norm\ closure); \ }$
\item[\ref{5.6}.4)] $i(\overline{co}(A\cup C)) =\max\{i(A),i(C)\}$ for all $A,C$ in $cwk(X)$;
\item[\ref{5.6}.5)] If \, $0\in{A}$, then $i(A)\geq 0$.
\end{description}
\end{thm}

    The embedding $i$  allows to reduce  the  Henstock (resp. McShane) integrability of multifunctions to  the    Henstock (resp. McShane) integrability of functions by embedding the family   $cwk(X)$  into the Banach space $C(\Omega)$.  Since
 $i(cwk(X))$ is a closed cone of $C(\Omega)$,   a multifunction $\vG:[0,1]\to cwk(X)$  is  Henstock or  variationally Henstock (resp. McShane or variationally McShane)
   integrable  if and only if the single valued function $i\circ \vG: [0,1] \rightarrow C(\Omega)$ is Henstock or  variationally Henstock (resp. McShane or variationally McShane) integrable in the usual sense.
The key point is that $i(cwk(X))$ is a closed cone. Consequently,  if $z\in{C(\Omega)}$ is the value of the  integral of $i\circ\vG$, then there exists a set $K\in{cwk(X)}$ with $i(K)=z$.\\

Observe that it follows directly from the definitions that if $i:cwk(X)\to {Y}$ is the  R{\aa}dstr\"{o}m embedding into a Banach space, then a multifunction $\vG:[0,1]\to{cwk(X)}$ is G-integrable if and only if $i(\vG)$ is G-integrable (G stands for any of the gauge integrals). In \cite[Proposition 2.6, Corollary 2.7]{CASCALES2} the authors proved it for Birkhoff.
Thanks to the embedding,   the following result, similar to that for single-valued function, is now obvious:

\begin{prop}\label{bvsp}
Let $\vG: [0,1] \rightarrow cwk(X)$ be a Birkhoff integrable multifunction. Then $\vG$ is McShane integrable.
\end{prop}
\begin{proof}
It is enough to observe that $\vG$ is McShane integrable if and only if $i(\vG)$ is McShane integrable.
\end{proof}

Finally  ${\mathcal{S}}_H(\vG)\;[{\mathcal{S}}_{MS}(\vG)\,,\,
   \mathcal{S}_P(\vG)\,,\,  {\mathcal{S}}_{HKP}(\vG)\,,{\mathcal{S}}_B(\vG)\,,\,  {\mathcal{S}}_{vH}(\vG)\,,\;{\mathcal{S}}_{vMS}(\vG)\;]$ denotes the family of all scalarly measurable selections of $\vG$
    that are Henstock [McShane, Pettis, Henstock-Kurzweil-Pettis, Birkhoff, variationally Henstock, variationally McShane] integrable.\\

A useful tool to study the integrability of a  single-valued function or of a multifunction is the variational  measure associated to the primitive.

    \begin{deff}\label{vma} \rm
Given a finitely additive interval multimeasure $\Phi:\mathcal{I} \rightarrow cwk(X)$, a gauge $\delta$ and a set $E \subset [0,1]$, we define
$$Var(\Phi, \delta,E)=\sup\sum_{j=1}^p\|\Phi(I_j)\|,$$
 where the supremum is taken over all the $\delta$-fine Perron partitions $\{(I_j,t_j)\}_{j=1}^p$ with  $t_j \in E$ for $j=1,\dots,p$.
The  set function
$$V_\Phi(E):=\inf_{\delta}\left\{Var(\Phi,\delta,E):\delta\ \text{is a gauge on }E\right\}$$
   is called \textit{the variational measure generated by} $\Phi$. Moreover, we say that $V_\Phi$ is \textit{absolutely continuous} with respect to $\lambda$  and we write $V_\Phi \ll \lambda$ if for every $E \in \mathcal{L}$ with  $\lambda(E)=0$ we have $V_\Phi(E)=0$.
\end{deff}

In order to deduce more information from the integral of a multifunction, we state the following simple result.
\begin{prop}\label{subint}
Let $F:[0,1]\to cwk(X)$ be any Henstock integrable mapping and  $f$ be a Henstock integrable selection of $F$. Then, for every interval $I\in\mcI$, one has
 $$(H)\int_I f dt\in (H)\int_I F dt\quad\mbox{and}\quad V_{\Phi_f}(I)\leq{V_{\Phi_F}(I)}\,.$$
\end{prop}
\begin{proof}
 If $x^*\in{X^*}$, then $x^*f\leq{s(x^*,F)}$. Hence $(H)\int_Ix^*f(t)\,dt\leq(H)\int_Is(x^*,F(t))\,dt$. The Hahn-Banach theorem yields the required result.
\end{proof}
Observe that the previous result holds also for McShane and Birkhoff integrable multifunctions: in these cases, moreover, the conclusion is valid not only for intervals, but also for arbitrary measurable subsets.\\

It is well known \cite{Ca} that a variationally Henstock integrable vector valued function is
Bochner measurable. Using the embedding of $cwk(X)$ into $C(\Omega)$ it is possible to prove a similar result for a multifunction.

\begin{prop}\label{p0}
 Let $\vG:[0,1]\to{cwk(X)}$ be variationally Henstock integrable. Then $\vG$ is
Bochner measurable.
\end{prop}
\begin{proof}
 Since  $\vG$ is variationally-Henstock integrable, then also
the single valued function
$\gamma:=i\circ \vG$ is so in $i(cwk(X))$.  By \cite[Theorem 9]{Ca} $i\circ \vG$ is
Bochner measurable.
\\
We may assume that $\gamma$ is separably valued. For each $n\in\N$ let
$\{y_{n,k}:k\in\N\}\subset \gamma([0,1])$ be a countable $1/n$-net of
$\gamma([0,1])$.
For each $n,k$  let
$$\gamma_n'(t)=y_{n,k} \,\, \mbox{if} \,\,\  t\in
\gamma^{-1}(B(y_{n,k},1/n)),$$
Then $(\gamma_n')_n$ is a sequence of countably valued functions that is
uniformly converging to $\gamma$.\\
 For each $n$ let $k_n$ be such that
$\lambda\left(\bigcup_{i=1}^{k_n}\gamma^{-1}(B(y_{n,i},1/n))\right)>1-1/n$.
We define
\[\gamma_n(t) = \left\{ \begin{array}{ll}
\gamma_n'(t) & t \in \bigcup_{i=1}^{k_n}\gamma^{-1}(B(y_{n,i},1/n))\\
& \\ \mbox{a fixed point} \, y_n\in{i(cwk(X))} & \mbox{otherwise}.
\end{array} \right. \]
Since the values of $\gamma_n$ belong to $\gamma([0,1])= i \circ
\Gamma ([0,1])$ and the embedding is one to one then each  function $\gamma_n$ defines in a unique way a simple
multifunction $\vG_n:[0,1]\to{cwk(X)}$ such that
$i\circ\vG_n(t)=\gamma_n(t)$ for every $t \in [0,1]$. Clearly
$d_H(\vG_n(t),\vG(t))=\|\gamma_n(t)-\gamma(t)\|\longrightarrow 0$ in
 $\lambda$-measure. Hence, there is a subsequence converging a.e. to $\vG$.
 In conclusion, $\vG$ is the a.e. limit of a sequence of simple multifunctions, and  therefore it is Bochner measurable.
\end{proof}

\section{Selections of $cwk(X)$-valued multifunctions}\label{five}

When $X$ is an arbitrary Banach space, then it is well known that each scalarly measurable selection of a Pettis (resp.  Henstock-Kurzweil-Pettis) integrable multifunction  $\vG:[0,1]\to{cwk(X)}$ is also Pettis (resp.  Henstock-Kurzweil-Pettis) integrable (see \cite{ckr,mu4,dm2}).\\ For the Henstock integral the behavior is different. If $X$ is separable and we consider $ck(X)$-valued multifunctions, every measurable selection of a Henstock (resp. McShane) integrable multifunction is Henstock (resp. McShane) integrable (see \cite{dm7}). This essentially depends on the fact that for functions taking values in a separable Banach space the Pettis and the McShane integrability coincide. \\ The answer is also affirmative if we consider  $cwk(X)$-valued multifunctions taking values in any  Banach space X with the property that the Pettis and the McShane integrability coincide, as the following proposition shows (see \cite{dp,r,apr} and the bibliography inside,  for the Banach spaces with such
 a property).

\begin{prop}\label{p1}
 Let $X$ be a  Banach space with the property that the Pettis and the McShane integrability coincide.  Then for any Henstock (resp. McShane) integrable multifunction $\vG:[0,1]\to{cwk(X)}$, every scalarly
measurable selection of $\vG$ is Henstock  (resp. McShane)  integrable.
\end{prop}
\begin{proof}   Let  $\vG:[0,1]\to{cwk(X)}$ be  Henstock integrable. According to \cite[Theorem 3.1]{dm} there exists $f\in{\mathcal S}_H(\vG)$ and so by \cite[Theorem 1]{dm2}  $\vG(t)=G(t)+f(t)$, where  $G:[0,1]\to cwk(X)$  is a Pettis integrable multifunction. Now let $h:[0,1]\to X$ be any scalarly measurable selection of $\vG$. So $h=g+f$, where g is a scalarly measurable selection of $G$. Since the multifunction $G$ is  Pettis integrable,  $g$ is also Pettis integrable (see \cite{ckr} or \cite{mu4}). Then by the hypothesis $g$ is also McShane (and then Henstock) integrable. This gives  that  $h$ as a sum of two Henstock integrable functions, is Henstock integrable. In case the multifunction $\vG$ is   McShane  integrable, then it is enough to consider the decomposition with a McShane integrable selection of  $\vG$ (again, by \cite[Theorem 3.1]{dm} one has ${\mathcal S}_{MS}(\vG)\neq\emp$).
\end{proof}

If $X$ is a general Banach space,  and the  multifunction is  $cwk(X)$ valued, the previous assertion is false as we are showing in the next proposition. To do it we  use an example given  in \cite[Theorem 3.7]{apr}  under ZFC,  of a scalarly negligible function,  which is not McShane integrable. We recall that a family $\mathcal{F}$ of finite subsets of $[0,1]$ is said to be {\it MC-filling} on $[0,1]$ if it is hereditary (i.e. if $G \in \mathcal{F}$ whenever $G\subset F \in \mathcal{F}$) and there exists $\ve >0$ such that for every countable family $(A_i)$ of disjoint sets whose union is $[0,1]$, there is $F\in \mathcal{F}$ such that
$\lambda^*\left(\cup \{A_i: \ F\cap A_i\neq \emptyset \} \right)>\ve,$
where $\lambda^*$ is the outer Lebesgue measure.

\vspace{2ex}
\begin{prop}\label{p2}
 There exist a reflexive Banach space $Y$ and a  variationally McShane integrable multifunction $\vG:[0,1]\to{cwk(Y)}$  such   that $0 \in \vG(t)$ for each $t \in [0,1]$ and $\vG$ possesses a scalarly measurable selection that   fails to be Henstock integrable.
\end{prop}
\begin{proof}
  Let  $\mathcal{F}$ be a compact MC-filling family on $[0,1]$, existing by  \cite[Theorem 3.5]{apr} and let $X$ be equal to $C(\mathcal{F})$ (the space of continuous functions on $\mathcal F$).  The function $f:[0,1]\to{C(\mathcal{F})}$ defined by
\[f(t)(F)=1_F(t), \,\, \mbox{for}\,\,
t \in [0,1] \,\, \mbox{and} \,\,  F \in \mathcal{F}\]
 is scalarly negligible and not McShane integrable. Moreover,  there exists a reflexive Banach  space $Y$, a one-to-one linear continuous mapping $T: Y\to{C(\mathcal{F})}$ and a function $g:[0,1]\to Y$ such that: $f([0,1]) \subset T(B(Y))$ and $T\circ g=f$. So also the function $g$ is scalarly negligible and not McShane (and then not Henstock, due to \cite{f1994a})  integrable and $g([0,1]) \subset B(Y)$.
 Since $Y$ is a reflexive Banach space, the unit ball $B(Y)$ is a convex weakly compact set of $Y$, not norm compact.
 Now let us consider the constant multifunction  $\vG: [0,1] \to{cwk(Y)}$   defined by
\[ \vG(t):= B(Y),\,\,  \mbox{for every}\,\, t \in [0,1].\]
 Clearly $\vG$ is variationally McShane
 integrable and  $g$  is a scalarly measurable selection of $\vG$  which fails to be Henstock integrable.
 \end{proof}

 \begin{prop}\label{p5}
   Each variationally Henstock  integrable multifunction $\vG:[0,1]\to{cwk(X)}$ possesses a
 strongly  measurable selection.
\end{prop}
\begin{proof}
 By Proposition \ref{p0} we know that $\vG$ is
Bochner measurable.
 We get the thesis applying \cite[ Theorem 2.9]{himmel} which states that: in a metric space any
Bochner measurable multifunction taking convex closed bounded values  has
strongly measurable selections. In \cite{himmel} the metric is supposed to be bounded but the result is available without this assumption (see  \cite[Remark 3.7]{hansell}).
\end{proof}

A similar selection theorem can be stated for Birkhoff integrable mappings.
\begin{thm}\label{birkhoff}
Let $\Gamma:[0,1] \to cwk(X)$ be any Birkhoff integrable multifunction. Then, there exists at least one Birkhoff integrable selection.
\end{thm}
\begin{proof}
Since $\gamma:=i\circ \Gamma$ is Birkhoff integrable,
 by \cite[Remark 1]{nara} the $\gamma$ is McShane integrable with a measurable gauge, namely
there exists an element $H\in cwk(X)$ such that, for each $\vp>0$ a measurable gauge $\delta$ can be found, such that, as soon as ${\mathcal P}=\{(t_j,I_j):j=1,...,n\}$ is any $\delta$-fine partition of $[0,1]$, it holds
$$ \| \sigma( \gamma, {\mathcal P})-i(H) \|_{\infty} \leq  \vp$$
Now, since  $\Gamma$ is McShane integrable, by \cite[Theorem 3.1]{dm}, it admits a McShane-integrable selection $g$. It only remains  to prove that, in the condition of  McShane integrability of $g$,  measurable gauges are involved. To this aim it is enough to repeat almost  verbatim the proof of  \cite[Theorem 3.1]{dm}, taking into account that, for each $x^*$ in the dual space, the real-valued map $x^*g$ is Lebesgue integrable, hence also Birkhoff integrable  again using \cite[Remark 1]{nara}.
\end{proof}

\begin{rem}\label{r1}
\rm
The example given in Proposition \ref{p2} shows that in general not all scalarly measurable selections of a Birkhoff integrable multifunction are Birkhoff integrable. \\
If $X$ is separable, then thanks to Pettis measurability Theorem,   the answer is contained in \cite[Proposition 3.1(ii)]{CASCALES2}. In fact in this case $\vG$ admits
strongly measurable selections each of them being Birkhoff integrable, and for every measurable set $A$ it is
\[ (B) \int_A \vG dt = \left\{ \int_A f dt,\hskip.4cm  f\in {\mathcal{S}}_B(\vG) \right\}.\]
In non-separable case we have that if $\vG:[0,1]\to{cwk(X)}$  is Birkhoff integrable,   then   $\vG$ is Pettis integrable by Proposition \ref{bvsp}.
 So each scalarly measurable selection of $\vG$ is Pettis integrable.
Birkhoff integrability can be obtained if the range of $\vG$ is separable, or if the selection is strongly measurable,
thanks to  \cite[Corollary 5.11]{pettis}.
\end{rem}

We recall  that a  multifunction $\vG:[0,1]\to{cwk(X)}$ is said to be {\it integrably bounded} if there is a scalar valued function $h\in{L_1[0,1]}$ such that $\|\vG(t)\|\leq|h(t)|$ for almost all $t\in[0,1]$.
 \begin{prop}\label{p4}
 Let $\vG:[0,1]\to{cwk(X)}$ be a scalarly measurable multifunction. Then the following conditions are equivalent:
 \begin{enumerate}
 \item[\bf (i)] $\vG$ is variationally McShane integrable;
  \item[\bf (ii)]  $i(\vG) \in  L_1(\lambda, Y)$;
  \item[\bf (iii)]  $\vG$ is  Bochner measurable  and integrably bounded.
  \end{enumerate}
 \end{prop}
\begin{proof}
 The equivalence of (i) and (ii)  follows from \cite{dm11}, where it has been proven that a Banach space valued
 function is variationally McShane integrable if and only if it is Bochner integrable. Consequently, $\vG$ is
 variationally McShane integrable if and only if  $i(\vG) \in  L_1(\lambda, Y)$. The equality
 $\|i \circ\vG(t)\|=\|\vG(t)\|$ and the fact that $i(\vG)$ is Bochner measurable if and only if $\vG$ is such,
 yield the equivalence of (ii) and (iii).
\end{proof}

  \begin{thm}\label{t10}
Each strongly measurable selection of a  variationally  McShane integrable multifunction $\Gamma:[0,1] \to cwk(X)$
  is vMS-integrable (= Bochner integrable).
  In particular, if  $X$ is separable, then  each scalarly measurable selection of $\vG$ is Bochner integrable.
  Consequently, $\vG$ is Aumann integrable and the integrals coincide.
\end{thm}
\begin{proof} That follows at once from Propositions \ref{p4} and \ref{p5} and the inequality
$\|f(t)\| \leq \|\vG (t)\|$ that is valid for any selection $f$ of $\vG$.
\end{proof}
 We would like to observe that in general  the Aumann integrability of a  measurable $cwk(X)$ valued
 multifunction does not imply its variational McShane integrability. In fact it is enough to consider
 Example \ref{ex1}: the  multifunction $G$   is not variationally McShane integrable, but is Aumann integrable,
 since the null function is a selection.
Moreover in \cite[Example 1]{bs2004} another example is given which is useful  for this purpose.
The example is reported here  for the sake of completeness and compared with the variational McShane integrability:
\begin{ex} \label{ex10}
\rm
Let $X= l^2(\mathbb{N})$; for every $A \subset \mathbb{N}$ we consider
\begin{eqnarray*}
U_A = \{ x \in X : \|x \| \leq 1, \mbox{~and~} x_n = 0 \mbox{~if ~} n
\not\in A \} = \{1_A x : \| x\| \leq 1\},
\end{eqnarray*}
where $(1_A x)_n = 1_A (n) x_n$. If $A \not=B$ then $d_H (U_A,U_B) \geq 1$ and
so the set $\{U_A, A \subset \mathbb{N}\}$ is not separable.\\
Let $\Omega = [0,1[$ and for every $\omega \in \Omega$ let $0,\omega_1
\cdots \omega_n \cdots$ be its dyadic representation, namely
$\omega_1 = 1$ ~iff~ $\omega \in [1/2,1[$, $\omega_2 = 1$ ~iff~ $\omega \in
[1/4,1/2[ \cup [3/4,1[$, etc.
We set $B_1= [1/2,1[, B_2 = [1/4,1/2[ ~\cup~ [3/4,1[$, etc.\\
Let $F(\omega) = U_{A(\omega)}$ where $A(\omega) = \{n \in
\mathbb{N} : \omega_n =1\}$. $F$ is integrably bounded, takes weakly
compact and convex values and its support function $s(y,F(\omega))$
is measurable since it is the limit of
simple functions. In particular all scalarly measurable selections are Bochner integrable.\\
From \cite[Proposition II.2.39]{hp} $F$ is Effros measurable, but
for every $\lambda$-null set $N,$ the set  $F(\Omega \setminus N)$
is not separable in the $d_H$-metric topology. Then immediately it
follows that $F$ is not Bochner  measurable and then it cannot be
variationally McShane integrable by Proposition \ref{p4}.
\end{ex}

 \begin{thm}\label{t1}
 If $\vG:[0,1]\to{ck(X)}$ is Bochner measurable and
Henstock  integrable ($\vG:[0,1]\to{cwk(X)}$ is variationally McShane  integrable),
then each scalarly measurable selection of $\vG$ is
 Henstock (variationally McShane) integrable.
Moreover in the second case
$\vG$ is Aumann integrable and the two integrals coincide.
\end{thm}
\begin{proof}
By \cite[Theorem 3.6]{dms} the range of $\vG$ is almost separably valued in $X$. Assume for simplicity that $\wt{X}$ is a  separable subspace of $X$ containing the range of $\vG$. It follows from \cite[Theorem 2.(iv)]{dm7} that all scalarly measurable selections of $\vG$ are Henstock integrable.
In case of variationally McShane integrable $\vG$, we apply Theorem \ref{t10} to  $\wt{X}$.
\end{proof}
\begin{prop}\label{pt1}
If $\vG:[0,1]\to{cwk(X)}$ is Birkhoff integrable and its range is almost separably valued then each scalarly measurable selection of
 $\vG$ is Birkhoff integrable, and then also McShane and  Henstock  integrable.
\end{prop}
\begin{proof}
Obviously $\vG$ is Pettis integrable and then, by \cite[Corollary 2.3]{ckr},  each scalarly measurable selection $f$ of $\vG$ is
Pettis integrable. Since the range of $\vG$ is almost separably valued the same occurs for $f$, so $f$ is McShane integrable.
Since $f$ is scalarly measurable and its range is separable, by \cite[Corollary 3]{nara}, it is Bochner measurable. But then $f$ satisfies the hypotheses of \cite[Theorem 7]{nara} and so it is Birkhoff integrable (see also the Remark \ref{r1}).
\end{proof}

The previous results, and in particular Theorem \ref{t10}, give raise to the following problem.
\begin{que} {\rm Does there exist at least one variationally Henstock   integrable selection of a $cwk(X)$-valued   variationally Henstock integrable multifunction?}
\end{que}

 If  $X$ is an arbitrary Banach space, we do not know the answer to the previous question. But for    Banach spaces
possessing the Radon-Nikod\'{y}m property   there exist variationally Henstock   integrable selections of  $cwk(X)$-valued   variationally Henstock integrable multifunctions, as it is stated in the next theorem.\\

  To present a proof, we need the notion of the variational  measure associated to the primitive
given in Definition \ref{vma}.\\

\begin{thm}\label{t11}
 Let $X$ be a Banach space with the Radon-Nikod\'{y}m property and let  $\vG:[0,1]\to {cwk(X)}$ be a variationally Henstock integrable multifunction. Then every strongly measurable selection of $\vG$ is variationally Henstock  integrable.
\end{thm}
\begin{proof}
Thanks to  Proposition \ref{p5} we know that there exist
 strongly measurable  selections of $\vG$.  Let $f$ be one of such selections. Since $\vG$ is variationally Henstock integrable, it is also Henstock-Kurzweil-Pettis integrable.  Therefore by \cite[Proposition 1.5]{dm5}  $f \in {\mathcal{S}}_{HKP}(\vG)$.\\ Let us denote by $F$ its
Henstock-Kurzweil-Pettis primitive and by $\vPh_{\vG}$ the Henstock primitive of $\vG$.
 By \cite[Proposition 3.3.1]{tesip}
we have $V_{\vPh_{\vG}} \ll \lambda$. Hence also $V_{F} \ll \lambda$. Since $X$ has
the Radon-Nikod\'{y}m property, by \cite[Theorem 3.6]{BDpM2} we infer that $F$ is the  variational Henstock primitive of $f$. Therefore  $f$ is variationally Henstock   integrable.
\end{proof}

 Comparing   the previous theorem  with Proposition \ref{p2}, we see that the condition of strong measurability of the selection cannot be relaxed.  Indeed,
in the example given in \ref{p2} $X$ is reflexive, and it is shown that a {\em constant} $cwk(X)$-valued multifunction has a scalarly measurable selection that is not even Henstock-integrable.\\

Another  partial answer    is the following:
\begin{thm}\label{selvm}
 If  $\vG:[0,1]\to {cwk(X)}$ is a bounded multifunction which
 is variationally Henstock and Pettis integrable, then ${\mathcal{S}}_{vMS}(\vG)\not=\emp$.
\end{thm}
\begin{proof}
Thanks to Proposition \ref{p5} the multifunction $\vG$ admits a
strongly measurable selection $f$ which is  Pettis integrable by \cite[Corollary 2.3]{ckr}. Let $\mu_f$ be its
indefinite Pettis integral. Then $R(\mu_f)$,   the range of $\mu_f$,   is relatively norm compact
(classical result) and the variation $|\mu_f|$ of $\mu_f$ is of $\sigma$-finite variation (cf. \cite[Theorem 4.1]{mu}).\\
Now, since $f$ is a selection of $\vG$ which  is bounded, then there exists $M > 0$
such that for every $E \in \mcL$ it is $|\mu_f| (E) \leq M \lambda(E)$, so
$\mu_f$ is moderated. Then  thanks to \cite[Lemma 2]{dm11} $f$ is  variationally McShane integrable.
\end{proof}

This result could be extended to multifunctions $F$ of the following type:
\begin{eqnarray*}
 \vG:[0,1]\to {cwk(X)} :   \exists \,\, M_n \in \mathbb{R} \mbox{\,\, and \,\,}  O_n =O_n^{\circ}  \mbox{\,\, such that\,\,\, }
 V_{\Gamma}(\cap_n O_n^c)=0\ \ {\rm and}\  \left\| \vG (t) \cdot 1_{O_n} (t) \right\| \leq M_n.
\end{eqnarray*}
These assumptions cannot be further weakened
since in \cite[Example 1]{mpot}
an example of a Birkhoff and strongly measurable function is given,
whose variational indefinite integral is not moderated on any open interval.
\\

It is well known that if $\vG$ is a $cwk(X)$ valued multifunction defined on a complete probability space (resp. $[0,1]$), and  all its scalarly measurable selections  are Pettis integrable (resp. Henstock-Kurzweil-Pettis integrable), then $\vG$ is Pettis integrable (resp. Henstock-Kurzweil-Pettis integrable) (\cite{ckr} or \cite{mu4} (resp. \cite{dm2})). A similar result is valid in case of a separable $X$ and a ck(X)-valued  Henstock  integrable  multifunction (\cite[Theorem 2]{dm7}).\\

The following proposition shows that the above assertion  remains true also for the McShane integral.

\begin{prop}\label{p3}
 Let $X$ be a separable Banach space and let  $\vG:[0,1]\to{ck(X)}$ be a  Bochner measurable  multifunction.
\begin{enumerate}
\item[\bf (i)]
 If all   measurable selections of  $\vG$ are McShane  integrable, then $\vG$ is McShane  integrable.
\item[\bf (ii)]
 If all   measurable selections of  $\vG$ are Birkhoff  integrable, then $\vG$ is Birkhoff  integrable.
\end{enumerate}
\end{prop}
\begin{proof}
 Since all  measurable selections of  $\vG$ are McShane  integrable, they are also Pettis and Henstock integrable. So by \cite[Theorem 4.2]{ckr} and \cite[Theorem 3.8]{ckr2}, $\vG$ is Pettis integrable, and by \cite[Theorem 2]{dm7}, $\vG$ is Henstock integrable. An application of \cite[Theorem 3.4]{dm} gives us the McShane  integrability of $\vG$.
The second statement  follows from  \ref{p3}.i) and \cite[Corollary 4.2]{BS2011}.
\end{proof}

 \begin{rem} {\rm We would like to point  out that even in the case of separable Banach spaces, the assertion of Proposition \ref{p3} is false if we consider  $cwk(X)$ valued multifunctions. In fact $\vG$ may fail  to be even  Henstock integrable. Indeed let $X$ be a separable Banach space without the Schur property. Then, following the proof  of \cite[Theorem 2.1]{ckr1} it is possible to construct a Pettis integrable multifunction $\vG :[0,1]\to{cwk(X)}$, such that $i\circ\vG$  is not scalarly measurable. Since $\vG$ is Pettis integrable, then (see \cite{ckr} or \cite{mu4}), each measurable selection of $\vG$ is Pettis, so  McShane integrable and then Henstock integrable. But $\vG$ cannot be Henstock integrable (hence neither McShane). In fact in such a case also $i\circ\vG$ would be Henstock integrable, whereas $i\circ\vG$ fails to be scalarly measurable.}
 \end{rem}

\section{Variationally Henstock and McShane integrability of $cwk(X)$-valued multifunctions}\label{three}

Some results concerning variationally Henstock and McShane integrable multifunctions are collected here.

\begin{prop}\label{prop1}
Let $G : [0, 1] \to cwk(X)$ be variationally Henstock integrable. If $0\in G(t)$ a.e.,
then $G$ is Birkhoff integrable.
\end{prop}
\begin{proof}
 Let $i$ be the embedding of $cwk(X)$ into $C(\Omega)$ of Theorem \ref{5.6}.
Then we just have to prove that $i(G)$ is Birkhoff integrable.\\
 To this aim, we observe that $i(G)$ is variationally Henstock integrable, hence Bochner-measurable, thanks to Proposition \ref{p0}. This implies that $i(G)$ is also Riemann-measurable, according with \cite[Theorem 1]{nara}. So, in order to prove that $i(G)$ is Birkhoff integrable, it is enough to show that it is McShane-integrable, thanks to \cite[Theorem 7]{nara}.
 Since $0\in G(t)$ one has that $i(G(t))$ is non-negative for almost all $t\in [0,1]$, and $i(G)$ is variationally Henstock integrable.\\
Then, thanks to \cite[Corollary 9 (iii)]{f1994a}, it will be sufficient to prove
convergence in $C(\Omega)$ of all series of the type $\sum_n (vH) \int_{I_n}i(G)$, where $(I_n)_n$ is any sequence of pairwise non-overlapping subintervals of $[0,1]$.\\
The map $\Psi(E):=(vH)\int_E G$ is defined, nonnegative and finitely additive on the algebra $\mathcal{H}$ generated by all intervals.
By \cite[Corollary 3.1]{BDpM2} $V_{\Psi} \ll \lambda.$
Since $0 \in G(t)$ then $s(x^*, \Psi (E)) \geq 0$ for every $x^* \in X^*$ and every $E \in \mathcal{H}$, and so by \cite[Theorem 4.6]{dpp}
 the map $\Psi$ can be extended to $\mathcal{L}$ in a $\sigma$-additive way; let $\widetilde {\Psi}: \mathcal{L} \to cwk(X)$  be its extension.
So, fixed   any sequence of pairwise non-overlapping subintervals $(I_n)_n$ of $[0,1]$, let $E = \cup_n I_n$.
Then
$i(\widetilde{\Psi})(E) := \sum_n(vH)\int_{I_n}i(G) \in C(\Omega).$
\end{proof}

\begin{thm}\label{t4}
Let $\vG:[0,1]\to cwk(X)$ be a  variationally Henstock integrable multifunction. If
${\mathcal{S}}_{vH}(\vG)\not=\emp$ (this is fulfilled in case of $X$ possessing RNP, by Proposition \ref{p5} and Theorem \ref{t11}), then for every
$f\in{\mathcal{S}}_{vH}(\vG)$  the multifunction $G:[0,1]\to cwk(X)$
defined by $\vG(t)=G(t)+f(t)$ is Birkhoff integrable;
\end{thm}

\begin{proof}
Let $f\in\mcS_{vH}(\vG)$ be fixed. Define  $G:[0,1]\to{cwk(X)}$ by $G(t):=\vG(t)-f(t)$. Then $G$ is also variationally Henstock integrable (in $cwk(X)$)  and   $0\in G(t)$ for every $t \in [0,1]$.  By Proposition  \ref{prop1} the multifunction $G$ is  Birkhoff integrable.
\end{proof}

The next two results generalize \cite[Theorem 3.4]{dm}, proved there for $cwk(X)$-valued multifunctions with compact valued integrals.
\begin{thm}\label{t3}
Let $\vG:[0,1]\to {cwk(X)}$ be a vH-integrable multifunction.
 If
${\mathcal{S}}_{vH}(\vG)\not=\emp$ (this is fulfilled in case of $X$ possessing RNP, by Proposition \ref{p5} and Theorem \ref{t11}), then    the following conditions are equivalent:
\begin{enumerate}
\item[ \it \bf (a)] $\mcS_{vH}(\vG)\subset\mcS_{MS}(\vG)$;
\item[ \it \bf (b)] $\mcS_{vH}(\vG)\subset\mcS_P(\vG)$;
\item[ \it \bf (c)] $\mcS_P(\vG)\not=\emp$;
\item[ \it \bf (d)] $\vG$ is   Pettis integrable.
\item[ \it \bf (e)] $\vG$ is McShane integrable.
\end{enumerate}
\end{thm}

\begin{proof}
 $\bs{(a)\Rightarrow(b)}$  is valid, because each McShane integrable
function is also Pettis integrable (\cite[Theorem 8]{f1994a}).
\\
$\bs{(b)\Rightarrow(c)}$ is obvious.
\\
 $\bs{(c)\Rightarrow(d)}$ Take  $f\in\mcS_P(\vG)$.  Since $\vG$ is
 Henstock integrable also, it is also HKP-integrable and so applying
\cite[Theorem 2]{dm2}  we obtain a representation $\vG=G+f$, where
$G:[0,1]\to{cwk(X)}$ is Pettis integrable in $cwk(X)$.  Consequently, $\vG$
is also Pettis integrable in $cwk(X)$ and so (d) is fulfilled.
\\
$\bs{(d)\Rightarrow (e)}$
In virtue of \cite[Theorem 3.1]{dm} $\vG$ has a
McShane integrable selection  $f$.    It follows from Theorem \ref{t4} that  the multifunction $G:[0,1]\to cwk(X)$
defined by $\vG(t)=G(t)+f(t)$ is McShane integrable.\\
$\bs{(e)\Rightarrow(a)}$
It is a consequence of Proposition \ref{p4} and \cite[Corollary 2.3]{ckr}, in fact,
if $f\in{\mathcal S}_{vH}(\vG)$, then
 $f\in {\mathcal S}_H(\vG)\cap{\mathcal S}_P(\vG)$ and so $f$ is McShane integrable.
\end{proof}

\begin{thm}\label{t2}
Let $\vG:[0,1]\to {cwk(X)}$ be an integrably bounded multifunction satisfying Theorem \ref{t3}.
 Then   all the statements  given in Theorem \ref{t3} are equivalent to
\begin{enumerate}
\item[ \it \bf (f)]
$\vG$ is  variationally McShane integrable.
\end{enumerate}
\end{thm}
\begin{proof}
The equivalence of the four first conditions can be proved as in Theorem \ref{t3}.

$\bs{(e)\Rightarrow(f)}$  Assume that $\vG$ is integrably bounded by $h\in{L_1[0,1]}$ and let $i$ be the  R{\aa}dstr\"{o}m  embedding of $cwk(X)$ into a Banach space $Z$. Since $\vG$ is variationally Henstock integrable, also $i(\vG)$ is variationally Henstock integrable. If $M_{\vG}$ is the indefinite Pettis integral of $\vG$, then it is a measure in the Hausdorff metric, hence also $i(M_{\vG})$ is countably additive in $Z$. Moreover, if $I\in{\mcI}$, then
\begin{equation}\label{e1}
\langle{z^*,i(M_{\vG}(I))}\rangle=(HK)\int_I\langle{z^*,i(\vG)(I)}\rangle\,d\lambda\qquad\mbox{for every}\;z^*\in{Z^*}.
\end{equation}
But $i(\vG)$ is integrably bounded and so every $\langle{z^*,i(\vG)(I)}\rangle$ is Lebesgue integrable. Consequently, both sides of the equality (\ref{e1}) may be extended to scalar measures on $\mcL$. By the integral boundedness assumption we have also $|i(M_{\vG})|(E)\leq \int_E|h|\,d\lambda$ for all $E\in\mcL$. Thus, $i(M_{\vG})$ is moderated (in fact finite). Then,  thanks to \cite[Lemma 2]{dm11}, $i(\vG)$ is  variationally McShane integrable. This proves the required result.

$\bs{(f)\Rightarrow(e)}$ is obvious.
\end{proof}
A comparison with the Birkhoff integrability is the following
\begin{prop}\label{abs}
If $\vG$ is Bochner measurable and abs(Bi)-integrable, then $\vG$ is vMS-integrable.
\end{prop}
\begin{proof}
Thanks to  Theorem \ref{5.6} the function $i \circ \vG$ is Bochner measurable and abs(Birkhoff)-integrable then, applying \cite[Corollary 1]{mpot},  $i \circ \vG$ is vMS-integrable.
\end{proof}

\noindent While for the converse implication only Birkhoff-integrability could be obtained, see  \cite[Corollary 3]{mpot}.
\begin{cor}
If $X$ has RNP and $\vG$ is variationally Henstock integrable, then every strongly measurable selection of $\vG$ is  variationally Henstock and Birkhoff integrable.  Moreover, $\vG$ turns out to be Birkhoff integrable too.
\end{cor}
\begin{proof}
The first part is an easy consequence of Theorem \ref{t11} and \cite[Corollary 2.3]{ckr}, since by \cite[Corollary 5.11]{pettis}, for finite measure spaces Birkhoff integrability of strongly measurable functions is equivalent to their Pettis integrability.  The final assertion about $\vG$ follows from Theorem \ref{t4} since both the selection $f$ and the translated mapping $G$ are Birkhoff integrable.
\end{proof}

 We would like to observe that in general the equivalent conditions  in  Theorem \ref{t3}  do not imply the variational  McShane integrability, even  in case of single valued functions (see  \cite{dmar}).
Moreover, also Proposition \ref{prop1} cannot be obtained, if we consider  variational McShane integrability instead of   Birkhoff integrability,  as the following example shows.
Observe first that  if $\vG$ satisfies assumptions of Proposition \ref{prop1}, then by Proposition \ref{bvsp} is Pettis integrable since it is McShane integrable.   The multifunction of the following Example \ref{ex1} satisfies Proposition \ref{prop1}, so it will be vH-integrable, Pettis integrable but not vMS-integrable.

 \begin{ex}\label{ex1}{\rm Assume that  $\sum_nx_n$ is unconditionally but not absolutely convergent in $X$. We assume that $\|x_n\|\leq 1$, for all $n\in\N$. Moreover, let $I_n=:(2^{-n},2^{-n+1}),\,n\in\N$. We define $f:[0,1]\to{X}$ by
$$
f(t)=\sum_{n=1}^{\infty}2^nx_n 1_{I_n}(t).
$$
It has been proven in \cite{dmar} that $f$ is variationally Henstock and Pettis integrable but is not variationally McShane integrable.
Consider now the multifunction $G:[0,1]\to{ck(X)}$ defined by $G(t)={\rm conv}\{f(t),0\}$.

{\bf Claim 1 } $G$ is variationally Henstock integrable.

\begin{proof}
For each $x^*$ we have $s(x^*,G(t))=(x^*f)^+(t)$ and so
if $I\in\mcI$, then
$$
\int_Is(x^*,G(t))\,dt=\int_I(x^*f)^+(t)\,dt=\int_I\left(\sum_{n=1}^{\infty}2^n\langle{x^*,x_n}\rangle^+1_{I_n}(t)\right)dt=
\sum_{n=1}^{\infty}2^n|I_n\cap{I}|\langle{x^*,x_n}\rangle^+
$$
where $\langle{x^*,x_n}\rangle^+=\langle{x^*,x_n}\rangle$ if $\langle{x^*,x_n}\rangle\geq 0$ and $0$ otherwise.
Let $\ve>0$ be fixed and $k\in\N$ be such that $\sup_{\|x^*\|\leq 1}\sum_{i=k}^{\infty}|\langle{x^*,x_i}\rangle|<\ve/2$.
We define a gauge setting
\[
\delta(t)= \left\{ \begin{array}{ll}
\min\{|t-2^{-n}|,|t-2^{-n+1}|\}&{\rm if} \ \  t\in{I_n}\\
&\\
\ve 2^{-2n-2}&{\rm if} \ \  t=2^{-n-1}  \\
&\\
\ve 2^{-k}& \ \ t=0
\end{array}\right.
\]
 Let $\{(J_1,t_1),\ldots,(J_p,t_p)\}$ be a $\delta$-fine Perron partition of $[0,1]$.
We have to evaluate the expression
\begin{eqnarray*}
&&  \sum_{i=1}^p\sup_{\|x^*\|\leq 1}\left|\int_{J_i}s(x^*,G(t))\,dt-s(x^*,G(t_i))|J_i|\right|\\
&=& \sum_{i=1}^p\sup_{\|x^*\|\leq 1}\left|\sum_{n=1}^{\infty}2^n|I_n\cap{J_i}|\langle{x^*,x_n}\rangle^+
-\sum_{n=1}^{\infty}2^n\langle{x^*,x_n}\rangle^+ 1_{I_n}(t_i)|J_i|\right|.
\end{eqnarray*}
If $t_i\in{I_{n_i}}$, then $J_i\subset{I_{n_i}}$ and we have
\begin{eqnarray*}
&& \left|\sum_{n=1}^{\infty}2^n|I_n\cap{J_i}|\langle{x^*,x_n}\rangle^+
-\sum_{n=1}^{\infty}2^n\langle{x^*,x_n}\rangle^+1_{I_n}(t_i)|J_i|\right|\\
&=& \left|2^{n_i}\langle{x^*,x_{n_i}}\rangle^+|J_i|-2^{n_i}\langle{x^*,x_{n_i}}\rangle^+|J_i|\right|=0.
 \end{eqnarray*}
If $t_i=2^{-n_i}$, then
\begin{eqnarray*}
&& \left|  \sum_{n=1}^{\infty}2^n|I_n\cap{J_i}|\langle{x^*,x_n}\rangle^+
-\sum_{n=1}^{\infty}2^n\langle{x^*,x_n}\rangle^+1_{I_n}(t_i)|J_i| \,\,\right|
= \left|\sum_{n=1}^{\infty}2^n|I_n\cap{J_i}|\langle{x^*,x_n}\rangle^+\right|\\
&&= 2^{n_i-1}|I_{n_i-1}\cap{J_i}|\langle{x^*,x_{n_{i-1}}}\rangle^+
+2^{n_i}|I_{n_i}\cap{J_i}|\langle{x^*,x_{n_i}}\rangle^+\\
&&\leq  2^{n_i-1}\ve 2^{2n_i-2}\langle{x^*,x_{n_{i-1}}}\rangle^+
+2^{n_i}\ve 2^{2n_i-2}\langle{x^*,x_{n_i}}\rangle^+\\
&&= \ve 2^{n_i-3}\langle{x^*,x_{n_{i-1}}}\rangle^++\ve 2^{n_i-2}\langle{x^*,x_{n_i}}\rangle^+<\ve 2^{-n_i}.
\end{eqnarray*}
If $t_i=0$, then
\begin{eqnarray*}
&& \left|\sum_{n=1}^{\infty}2^n|I_n\cap{J_i}|\langle{x^*,x_n}\rangle^+
-\sum_{n=1}^{\infty}2^n\langle{x^*,x_n}\rangle^+1_{I_n}(t_i)|J_i|\right|\\
&=&\sum_{n=1}^{\infty}2^n|I_n\cap{J_i}|\langle{x^*,x_n}\rangle^+\leq \sum_{n=k}^{\infty}2^n|I_n|\langle{x^*,x_n}\rangle^+\leq \sum_{n=k}^{\infty}|\langle{x^*,x_n}\rangle|\leq \ve/{2^{-k}}.
\end{eqnarray*}
So we can conclude that
$$\sum_{i=1}^p\sup_{\|x^*\|\leq 1}\left|\int_{J_i}s(x^*,G(t))\,dt-s(x^*,G(t_i))|J_i|\right|<2\vp.$$
\end{proof}
{\bf Claim 2 } $G$ is not variationally McShane integrable  and not abs-Birkhoff.
\begin{proof}
  Let $\delta:[0,1]\to (0,\infty)$ be an arbitrary gauge. Let $m\in\N$ be the smallest number satisfying the inequality $2^{-m}\leq\delta(0)$. For each $p>m$ we take a particular $\delta$-fine partition of $[0,1]$: 
  $$\mcP:=\{([0,2^{-p}],0), (I_p,0),\ldots,(I_{m+1},0), (J_1,s_1),\ldots,(J_q,s_q)\},$$ where we require each $(J_i,s_i)$
   to be only $\delta(s_i)$-small.
We have then
\begin{eqnarray*}
&& \sum_{(J,t)\in\mcP}\sup_{\|x^*\|\leq 1}\left|\sum_{n=1}^{\infty}2^n|I_n\cap{J}|\langle{x^*,x_n}\rangle^+
-\sum_{n=1}^{\infty}2^n\langle{x^*,x_n}\rangle^+1_{I_n}(t)|J|\right|\\
&\geq& \sum_{i=m+1}^p\sup_{\|x^*\|\leq 1}\left|\sum_{n=1}^{\infty}2^n|I_n\cap{I_i}|\langle{x^*,x_n}\rangle^+\right|
\\&=&
\sum_{i=m+1}^p\sup_{\|x^*\|\leq 1}2^i|I_i|\langle{x^*,x_i}\rangle^+=\sum_{i=m+1}^p\|x_i\|\longrightarrow\infty
\end{eqnarray*}
when $p\to\infty$.
 For the last part, since $G$ is Bochner measurable if it were abs-Birkhoff, then by Proposition \ref{abs} then $G$ would be variationally McShane-integrable.
\end{proof}}
\end{ex}

\section*{Acknowledgement}
The  paper was written while the third named author was visiting the
 Departments of Mathematics and Computer Sciences of the  University  of Palermo (Italy) thanks to the Grant  Prot. N. U 2014-001415 of GNAMPA -- INDAM (Italy).

The first and the last
 authors have been supported by University of Perugia -- Department of Mathematics and Computer Sciences -- Grant Nr 2010.011.0403 and, respectively, by Prin "Metodi logici per il trattamento dell'informazione" and    "Descartes".
The second author was supported by F.F.R. 2013 - University of Palermo (Italy).\\

\noindent This is a preprint version of an article published in Journal of Mathematical Analysis and Applications. The final authenticated version is available online at:\\
https://www.sciencedirect.com/science/article/abs/pii/S0022247X16300385

\end{document}